\documentclass[12pt]{article}

\usepackage{amsmath,amsfonts,amsthm, amssymb}

\newcommand{\red}[1]{{\color{red}{#1}}}
\newcommand{\blue}[1]{{\color{blue}{#1}}}

\usepackage{graphicx}
\usepackage{amsthm}
\usepackage{mathrsfs}
\usepackage{pgf,tikz}
\usepackage{cancel}
\usepackage{newclude}
\usepackage{multirow}

\def \Z {\mathbb Z}

\def \Cay {\mathrm{Cay}}

\newcommand{\comment}[1]{}

\newtheorem{theorem}{Theorem}[section]

\newtheorem{definition}[theorem]{Definition}

\newtheorem{lemma}[theorem]{Lemma}
\newtheorem{corollary}[theorem]{Corollary}
\newtheorem{conjecture}[theorem]{Conjecture}

\theoremstyle{definition}
\newtheorem{remark}[theorem]{Remark}
\newtheorem{example}[theorem]{Example}

\begin{document}

\title{On equitably 2-colourable  odd cycle decompositions}

\author{Andrea Burgess \thanks{andrea.burgess@unb.ca}\\
Department of Mathematics and Statistics,\\
University of New Brunswick \\Saint John, Canada \and
Francesca Merola  \thanks{francesca.merola@uniroma3.it} \\
Dipartimento di Matematica e Fisica,\\ Universit\`a Roma Tre,\\Rome, Italy }

\date{\today}
\maketitle
\begin{abstract}
    An $\ell$-cycle decomposition of $K_v$ is said to be \emph{equitably $2$-colourable} if there is a $2$-vertex-colouring of $K_v$ such that each colour is represented (approximately) an equal number of times on each cycle: more precisely, we ask that in each cycle $C$ of the decomposition, each colour appears on $\lfloor \ell/2 \rfloor$ or $\lceil \ell/2 \rceil$ of the vertices of $C$.
    In this paper we study the existence of equitably 2-colourable $\ell$-cycle decompositions of $K_v$, where $\ell$ is odd, and prove the existence  of such a decomposition for $v \equiv 1, \ell$ (mod $2\ell$).

\end{abstract}

\section{Introduction}

In this paper, we consider the existence of equitably 2-colourable $\ell$-cycle decompositions of $K_v$, where $\ell$ is odd.  We say that a graph $\Gamma$ {\em decomposes} into subgraphs $B_1, B_2, \ldots, B_t$ if the edge sets of the $B_i$ partition the edges of $\Gamma$.  In the case that $B_1 \simeq B_2 \simeq \cdots \simeq B_t \simeq B$, we speak of a {\em $B$-decomposition} of $\Gamma$.  In particular, an $\ell$-cycle decomposition of $\Gamma$ is a $B$-decomposition of $\Gamma$ where $B$ is an $\ell$-cycle $C_{\ell}$.

The study of cycle decompositions originates in the 19th century, with the classical results of Kirkman and Walecki on $3$-cycle and Hamilton cycle decompositions of the complete graph~\cite{Kirkman, Walecki}.  If $\Gamma$ decomposes into cycles, then each vertex of $\Gamma$ must have even degree. Thus, the complete graph $K_v$ admits a cycle decomposition only when $v$ is odd; for $v$ even, it is common to instead decompose the cocktail party graph $K_v-I$, which is formed from $K_v$ by removing the edges of a 1-factor $I$.  Where necessary, we will use the notation $K_v^*$ to denote $K_v$ if $v$ is odd and $K_v-I$ if $v$ is even.  
The existence of $\ell$-cycle decompositions of $K_v^*$ was solved by Alspach, Gavlas and \v{S}ajna~\cite{AlspachGavlas, Sajna}; see also~\cite{Buratti2003} for an alternative proof in the odd-cycle case.
\begin{theorem}[\cite{AlspachGavlas, Sajna}]
\label{AGS}
There exists a $C_{\ell}$-decomposition of $K_v^*$ if and only if $3 \leq \ell \leq v$ and $\ell$ divides $v\left\lfloor \frac{v-1}{2}\right\rfloor$.
\end{theorem}
\noindent
Note that in particular, there exists an $\ell$-cycle decomposition of $K_v$ whenever $v \equiv 1 \pmod{2\ell}$ or $\ell$ is odd and $v \equiv \ell \pmod{2\ell}$.

Decompositions of $K_v^*$ into cycles of varying lengths have also been studied.  In 1981, Alspach~\cite{Alspach1981} conjectured that the obvious necessary conditions for decomposing $K_v^*$ into cycles of given lengths are also sufficient.  This was finally proven in 2014 by Bryant, Horsley and Pettersson.
\begin{theorem}[\cite{BHP}]\label{BHPTheorem}
Let $n\geq 3$ be an integer and $\ell_1, \ell_2, \ldots, \ell_t$ be a list of integers (not necessarily distinct) with $3 \leq \ell_i \leq v$ for all $i \in \{1, \ldots, t\}$.  There is a decomposition of $K_v^*$ into cycles of lengths $\ell_1, \ell_2, \ldots, \ell_t$ if and only if $\ell_1+\ell_2+\cdots+\ell_t=v\left\lfloor\frac{(v-1)}{2}\right\rfloor$. 
\end{theorem}

The focus of this paper is a problem related to vertex colourings of cycle decompositions.  If $\mathcal{C}$ is a cycle decomposition of $\Gamma$, then a {\em $c$-colouring} of $\mathcal{C}$ is a function $\phi: V(\Gamma) \rightarrow S$, where $|S|=c$.  Informally, a $c$-colouring may be thought of as an assignment of $c$ colours to the vertices of $\Gamma$.  We will mainly consider the case $c=2$, and will denote the colours by {\em red} and {\em blue}.  

Of course, we generally wish such a colouring to satisfy additional properties concerning the cycles of $\mathcal{C}$.  For example, if we require that each cycle have at least two vertices coloured differently, the colouring is called weak; see~\cite{BurgessPike1, BurgessPike2, HorsleyPike} for results on weak colourings of cycle systems.  In this paper, we will consider equitable colourings.  We say that a $c$-colouring of an $\ell$-cycle decomposition $\mathcal{C}$ is {\em equitable} if each cycle contains either $\lfloor \ell/c\rfloor$ or $\lceil \ell/c \rceil$ vertices of each colour, i.e.\ the colours are as equally distributed as possible amongst the vertices of each cycle.  Note that the term {\em equitable colouring} also occurs in the literature in other contexts, in particular to refer to colourings in which the sizes of the colour classes differ by at most 1 (see, for instance,~\cite{TripleSystems}); our use of the word {\em equitable} follows the definition introduced in the Ph.D. thesis of Waterhouse~\cite{WaterhousePhD} (see also~\cite{ABW2007,ABLW2004,LefevreWaterhouse2005,Waterhouse2006}).

Equitable colourings of cycle systems of the complete and cocktail party graphs were considered in~\cite{ABW2007} and~\cite{ABLW2004}, while the papers~\cite{LefevreWaterhouse2005} and~\cite{Waterhouse2006} consider equitable colourings of complete multipartite graphs.  The main results of these papers restrict their attention to the case that the number of colours is 2 or 3, and the cycle length is small (at most 6).  In particular, in~\cite{ABW2007} and~\cite{ABLW2004}, the authors completely determine the existence of equitably $2$- or $3$-colourable $\ell$-cycle systems of $K_v$ and $K_v-I$ when the cycle length $\ell$ is $4$, $5$ or $6$. In~\cite{BM2021}, the present authors consider equitable 2-colourings of even cycle systems of $K_v-I$. 

In this paper, we consider equitable $2$-colourings of $\ell$-cycle systems of the complete graph $K_v$.  A straightforward counting argument (see~\cite{ABW2007, BM2021}) shows that no equitably $2$-colourable $\ell$-cycle decomposition of $K_v$ can exist if $\ell$ is even; hence necessarily we consider only odd cycle length $\ell$.  We focus on the case where $v \equiv 1, \ell$ (mod $2\ell$), in which there exists an $\ell$-cycle system of $K_v$ for any $\ell \geq 3$.  Our main result is the following.
\begin{theorem}\label{th:main}
Let $\ell>5$ be odd, and $v \equiv 1$ or $\ell \pmod{2\ell}$.  There is an equitably 2-colourable $\ell$-cycle decomposition of $K_v$.
\end{theorem}

In the case that the cycle length $\ell=p$ is prime, Theorem~\ref{AGS} asserts that a $p$-cycle system of $K_v$ exists if and only if $v \equiv 1$ or $p \pmod{2p}$.  Thus, Theorem~\ref{th:main}, combined with the results of~\cite{ABW2007} for cycle length $5$ and the nonexistence of an equitably $2$-colourable Steiner triple system~\cite{Rosa} of order greater than 3, gives a complete existence result for equitably 2-colourable $p$-cycle systems of $K_v$.
\begin{corollary}
Let $p$ be an odd prime.  There exists an equitably 2-colourable $p$-cycle decomposition of $K_v$ if and only if $p \geq 5$ and $v \equiv 1$ or $p\pmod{2p}$, or $p=v=3$.
\end{corollary}

Theorem~\ref{th:main} will be proved in Theorem \ref{th:1mod2l} for $v \equiv 1$ (mod $2\ell$), and in Theorem \ref{th:lmod2l} for $v \equiv \ell$ (mod $2\ell$).  As we will see, the proof of  Theorem \ref{th:1mod2l} relies on the existence of an equitably $2$-colourable $\ell$-cycle decomposition of $K_{2\ell+1}$ with colour classes of sizes $\ell$ and $\ell+1$, and this will be established in Section \ref{se:2ell+1}.  Note that, moreover, the existence of an equitably $2$-colourable $\ell$-cycle decomposition of $K_{4\ell+1}$ cannot be established via the methods used in proving Theorem \ref{th:1mod2l}, but requires a somewhat long ad hoc construction, which will be presented in Section  \ref{se:4ell+1}.

\section{The $2\ell+1$ and $4\ell+1$ case}

In this section, we prove the existence of an equitably $2$-colourable $\ell$-cycle decomposition of $K_{2\ell+1}$ (in Section~\ref{se:2ell+1}) and  $K_{4\ell+1}$ (in Section~\ref{se:4ell+1}).

The constructions in Sections \ref{se:2ell+1} and \ref{se:4ell+1} will use $2$-rotational cycle systems (see for instance \cite{Vietri}), i.e. cycle systems of order $2n+1$ admitting an automorphism having exactly
one fixed point and two cycles of length $n$ each.

The set of vertices will be identified with $(\Z_n\times \{0,1\}) \cup \{\infty\}$, with $\infty$ being the fixed point as usual in this situation, and we will often denote the vertex $(a,i), i=0,1$, with $a_i$. 
The colourings used will differ in the two cases, but share the property that vertices $a_0$ and $a_1$ will have different colours. The vertex $\infty$ will have colour blue.

A cycle system will be built via difference methods using a set of base cycles, which will be developed modulo $(n, -)$ 
i.e.\ $\pmod{n}$ with respect to the first coordinate (with the usual understanding that $\infty+i=\infty$); 
each base cycle $C$ used in the decompositions that follow will in this way give rise to an orbit of $d$ cycles, $d \mid n$. 
As usual, the difference arising from an edge of a cycle $a_i\sim b_j$  will be called {\em mixed} if $i\ne j$ and {\em pure} (or {\em $i$-pure}) if $i=j$. In particular, our convention is that the edge $a_0\sim b_1$ will give the mixed difference $b-a \pmod{n}$, and the  edge $a_i\sim b_i$ will give the $i$-pure difference $\pm(b-a) \pmod{n}$.

We point out that the constructions used here are similar to those used in \cite{WangCao, WangLuCao} to obtain almost resolvable odd $\ell$-cycle systems of order $2\ell+1$ and $4\ell+1$. In particular, a suitable colouring of the cycle systems built in \cite{WangLuCao}   gives an equitably $2$-colourable $\ell$-cycle decomposition of $K_{2\ell+1}$ when $\ell\equiv 1 \pmod{4}$.

\subsection{Case $v = 2\ell+1$}\label{se:2ell+1}

In this section we show the existence of an equitably $2$-colourable $\ell$-cycle decomposition of $K_{2\ell+1}$ with colour classes of sizes $\ell$ and $\ell+1$. As mentioned before, this cycle system will be 2-rotational and
will arise from a set of three base cycles to be developed 
modulo $(\ell, -)$, i.e. 
$\pmod{\ell}$ with respect to the first coordinate; 
the first cycle contains $\infty$ and 
has an orbit of length $\ell$, the second base cycle will have trivial stabilizer and give a length $\ell$ orbit, and the third cycle is stabilized by $\Z_\ell$ and has an orbit of length one.

To build the  cycle containing $\infty$ we will use a result on graceful labellings.  A {\em graceful labelling} of a path $P$ of length $k$ is an injective assignment of the integers $0, \ldots, k$ to the vertices of $P$ with the property that each integer $1, \ldots, k$ appears as the difference of the endvertices of an edge (see \cite{Gallian}).  For example, the ordered labelling $0,5,1,4,2,3$ gracefully labels a path of length $5$.

\begin{lemma}\label{le:graceful}
Consider a path with $h\ge 2$ vertices having labels $0,1,\dots,h-1$. There is a graceful labelling of this path, say $P_0$, such that one leaf has label $0$, and one, say $P_1$, in which one leaf has label $1$.
\end{lemma}

The proof follows, for instance, from Theorem 3 in \cite{Cat}.

\begin{theorem} \label{th:2l+1}
There is an equitably $2$-colourable $\ell$-cycle decomposition of $K_{2\ell+1}$ with colour classes of sizes $\ell$ and $\ell+1$.
\end{theorem}

\begin{proof}
Let $V(K_{2\ell+1}) = (\mathbb{Z}_\ell \times \{0,1\}) \cup \{\infty\}$; as already mentioned, we will often denote the vertex $(a,i), i=0,1$ with $a_i$ for brevity. We colour the vertices by using red on the vertices with an even label, and blue otherwise in the part labelled $0$, and reverse this colouring in the part with label $1$.  We colour $\infty$  blue. We will build a cycle system developing modulo $(\ell,-)$.

We can use Lemma \ref{le:graceful}   to construct
 an $\ell$-cycle $C_\infty$ having $\ell-1$ edges of $i$-pure differences $\pm 1,\pm 2,\dots, \pm((\ell-3)/2)$, $i=0,1$, and one of mixed difference 1, which is equitably $2$-coloured with respect to the colouring given above, as follows.

Set $h=(\ell-1)/2$ and apply Lemma \ref{le:graceful} to build paths $P_0$ on the $h=(\ell-1)/2$ vertices $0_0, \dots,(h-1)_0$, with the last leaf labelled $0_0$, and $P_1$ on $0_1, \dots,(h-1)_1$, with the first leaf labelled $1_1$. We obtain the cycle $C_\infty$  by concatenating $\infty, P_0$ and $P_1$.

Now consider the $\ell$-cycle

\[
C=\left(0_1, (\ell-2)_0, 1_1,(\ell-3)_0,2_1, \dots,(\frac{\ell+1}{2})_0,(\frac{\ell-3}{2})_1,(\frac{\ell-1}{2})_0,(\frac{\ell-1}{2})_1\right).
\]
The cycle $C$ contains an edge of each mixed difference except for $1$  and one edge of 1-pure difference $\pm(\ell-1)/2$.

Note that the $\ell$-set of cycles we obtain by developing this cycle  on first coordinates consists of equitably coloured cycles. 
Indeed the last two vertices will have different colours in $C$ and all its translates; the remaining vertices  form a run of $(\ell-1)/2$ consecutive 
$\pmod{\ell}$ $1$-vertices,  and $(\ell-3)/2$ consecutive 
$\pmod{\ell}$ $0$-vertices, starting from, in the cycle $C$, $0_1,\dots,(\ell-3)/2_1$ and $(\ell+1)/2_0, \dots,(\ell-2)_0$.
It can be checked that this pattern remains equitably coloured in all translates: this is fairly obvious for $C+i, i\ne 2,(\ell+3)/2$ and can be checked explicitly in these cases.

Now $C$ and $C_\infty$ between them contain an edge of each mixed difference and each pure difference except for the 0-pure difference $\pm(\ell-1)/2$,
and  $(\ell-1)/2$ is a unit in $\Z_\ell$; hence this last difference is covered by the cycle 
$$C_{(\ell-1)/2}=\left(0_0,((\ell-1)/2)_0,(2(\ell-1)/2)_0,(3(\ell-1)/2)_0,\dots\right).$$
\end{proof}

\begin{example}
For $\ell=9$, we can take the three cycles 
$$
\begin{array}{lll}
C_\infty & = & ({\color{blue}\infty}, {\color{red}2_0},{\color{blue}1_0},{\color{blue}3_0},{\color{red}0_0}, {\color{red}1_1},{\color{blue}2_1}, {\color{blue}0_1},{\color{red}3_1}),\\ 
C & = & ({\color{blue}0_1}, {\color{blue}7_0}, {\color{red} 1_1}, {\color{red} 6_0}, {\color{blue}2_1}, {\color{blue}5_0}, {\color{red}3_1}, {\color{red}4_0}, {\color{blue}4_1}),\\
C_4 & = & ({\color{red}0_0}, {\color{red}4_0}, {\color{red}8_0}, {\color{blue}3_0}, {\color{blue}7_0}, {\color{red}2_0}, {\color{red}6_0}, {\color{blue}1_0}, {\color{blue}5_0}).
\end{array}
$$

\end{example}

\subsection{Case $v=4\ell+1$}\label{se:4ell+1}

The aim of this section is to show an explicit construction for an equitably 2-colourable $\ell$-cycle decomposition of $K_{4\ell +1}$. 
We will prove the following.
\begin{theorem} \label{th:4l+1}
There is an equitably $2$-colourable $\ell$-cycle decomposition of $K_{4\ell+1}$ with colour classes of sizes $2\ell$ and $2\ell+1$.
\end{theorem}

As above, we will construct a $2$-rotational cycle system, this time developed modulo $(2\ell,-)$ from a set of six base cycles.
On part $0$ we will colour the vertices from $0$ to $\ell-1$ red and those from $\ell$ to $2\ell-1$ blue, and reverse this colouring in part 1. The vertex $\infty$ has colour blue. 

To help the reader, we will first present an example to clarify the general construction.

\begin{example}\label{ex:general4l+1}
Let $\ell=17$ and $v=69$: 
we identify the vertices of $K_{69}$ with $\Z_{34}\times\{0,1\}$, and denote the vertex $(a,i), i=0,1$ with $a_i$.

The cycle system has six base cycles which will be developed $\pmod{34}$. The first three cycles $C_p, C_0$ and $C_1$, have orbits of length $34$ and, in our example, are as follows.
\[
\begin{array}{ll}
C_p &=(\red{0_0},\red{1_0},\blue{17_0},\red{2_0},\red{8_0},\red{3_0},\red{7_0},\red{4_0},\blue{4_1},\blue{7_1},\blue{3_1},\blue{8_1},\blue{2_1},\red{10_1},\blue{0_1},\blue{1_1},\red{17_1})\\
C_0 &= (\red{6_0},\blue{9_1},\red{4_0},\blue{11_1},\red{2_0},\blue{13_1},\red{0_0},
\blue{6_1},\red{9_0},\blue{4_1},\red{11_0},\blue{2_1},\red{13_0},\blue{0_1},
\red{16_0},\blue{23_0},\red{22_1})\\
C_1 &= (\blue{6_1},\red{8_0},\blue{4_1},\red{12_0},\blue{2_1},\red{14_0},\blue{0_1},
\red{6_0},\blue{8_1},\red{4_0},\blue{12_1},\red{2_0},\blue{14_1},\red{0_0},
\blue{15_1},\red{22_1},\blue{21_0})\\
\end{array}
\]
Note that the differences from $C_0$ and $C_1$ cover all the mixed differences except for $0$ and $\ell=17$, and in $C_0$ and $C_1$ we also have one edge giving a pure difference, in the example the $0$-pure difference $\pm 7$ in $C_0$ and the $1$-pure difference $\pm 7$ in $C_1$. The mixed differences $0$ and $\ell=17$ appear in $C_p$, and the remaining differences from $C_p$ are all pure, namely $0$-pure differences 
$\pm\{1,3,4,5,6,15,16\}$
and $1$-pure differences 
$\pm\{1,3,4,5,6,8,10,16\}$.

Two base cycles $C_\infty^0, C_\infty^1$
 will contain the vertex $\infty$  and all remaining $0$-pure (resp. 1-pure) differences, except for $0$-pure difference $\pm 2$. 
These will be built using Lemma \ref{CyclesThroughInfinity}, which states that it is possible to build cycles $C_\infty^i, i=0,1$, that will contain $\infty$, have edges of $i$-pure difference $\ell$ and any $(\ell-3)/2$ other $i$-pure differences; these cycles will have orbit-length $\ell$.
In our case we have:

\[
\begin{array}{ll}
C_\infty^0 &=(\blue{\infty},\red{0_0},\red{14_0},\red{1_0},\red{13_0},\red{2_0},\red{12_0},\red{3_0},\red{11_0},\blue{28_0},\blue{20_0},\blue{29_0},\blue{19_0},\blue{30_0},\blue{18_0},\blue{31_0},\blue{17_0}),\\
C_\infty^1 &=(\blue{\infty},\blue{0_1},\blue{15_1},\blue{1_1},\blue{14_1},\blue{2_1},\blue{13_1},\blue{4_1},\blue{6_1},\red{23_1},\red{21_1},\red{30_1},\red{19_1},\red{31_1},\red{18_1},\red{32_1},\red{17_1}).
\end{array}
\]

We are left with the $0$-pure difference $\pm 2$ that gets covered by the cycle
$$C_2^0=(\red{0_0},\red{2_0},\red{4_0},\red{6_0},\red{8_0},\red{10_0},\red{12_0},\red{14_0},\red{16_0},\blue{18_0},\blue{20_0},\blue{22_0},\blue{24_0},\blue{26_0},\blue{28_0},\blue{30_0},\blue{32_0})$$ having orbit length 2.
\end{example}

\begin{proof}[Proof of Theorem~\ref{th:4l+1}]
    Let us now describe the general construction needed to prove the theorem.
As already mentioned, we will use a set of 6 base cycles, namely $\{C_0,C_1,C_p,C_\infty^0,C_\infty^1,C_2^0\}$, developed modulo $(2\ell,-)$.
The cycles $C_0$ and $C_1$ will cover all the mixed differences except for $0$ and $\ell$, and $C_i$ will cover a single $i$-pure difference $m\in \pm \{(\ell-5)/2,(\ell-3)/2,(\ell-1)/2,(\ell+1)/2\}$, $i=0,1$.  The construction of $C_0$ and $C_1$ is somewhat labourious, and it is described in Lemma \ref{le:C01}.

Let $S^0=\pm\{d_1,\dots,d_{(\ell-3)/2}\}$ be a set $0$-pure differences of size $(\ell-3)$ not containing $\pm 2,\pm m$ or $\ell$, and $S^1=\pm\{d_1,\dots,d_{(\ell-1)/2}\}$ a set $1$-pure differences of size $(\ell-1)$ not containing $\pm m$ or $\ell$.
The cycle $C_p$ will have on its edges the mixed differences $0$ and $\ell$, the $0$-pure differences in $S^0$, and the $1$-pure differences in $S^1$. Its construction is to be found in Lemma \ref{le:Cp}, and in Lemma \ref{le:l=3mod8} for the case $\ell\equiv3\pmod{8}$.

The cycles $C_0,C_1,C_p$ are all equitably coloured in the colouring described above.

The remaining pure differences will be found in the last three cycles.
The cycles $C_\infty^0,C_\infty^1$ are built applying Lemma \ref{CyclesThroughInfinity}, which gives us the freedom to use any $\ell-3$ $0$-pure and $1$-pure difference not already covered, and their development will also contain all edges through $\infty$. It is easily checked that these cycles and their translates are equitably coloured. 

The $0$-pure difference $\pm 2$ has not yet been used, and this difference can be found on the last base cycle $C_2^0=(0_0,2_0,4_0,\dots,(2\ell-2)_0)$, which is equitably coloured by construction, together with its translate. This completes the proof.  A summary of the differences appearing in each cycle appears in Table~\ref{DiffTable}.
\begin{center}
\begin{table}[ht]
\centering
\begin{tabular}{ |c|c|c| } 
 \hline
\multirow{2}{*}{mixed differences} & $0,\ell $ & $C_p$ \\ 
  \cline{2-3}
  &$ \{0,1,\dots2\ell-1\}\setminus\{0,\ell\}$ & $C_0$ and $C_1$\\
  \hline
 \multirow{4}{*}{$0$-pure differences} & $\pm m $ & $C_0$ \\ 
  \cline{2-3}
  & $\pm 2$ & $C_2^0$ \\ 
  \cline{2-3}
  & $S^0$ &$C_p$\\
  \cline{2-3}
  & $\Z_{2\ell}\setminus(S^0\cup\{0,\pm 2,\pm m\})$ & $C^0_\infty$\\
 \hline

\multirow{3}{*}{$1$-pure differences} & $\pm m $ & $C_1$ \\ 
\cline{2-3}
  &  $S^1$  &$C_p$\\
  \cline{2-3}
  &  $\Z_{2\ell}\setminus(S^1\cup\{0,\pm m\})$ & $C^1_\infty$\\
  
  \hline
  \hline
\multirow{2}{*}{edges containing $\infty$} &  & $C_{\infty}^0$ and $C_{\infty}^1$\\ & &and their translates \\ 
  
  \hline
\end{tabular}
\caption{The differences used in $C_0, C_1, C_p, C_2^0 C_{\infty}^0$ and $C_{\infty}^1$.\label{DiffTable}}
\end{table}
\end{center}

\end{proof}

Let us show firstly how to build two base cycles $C_0$ and $C_1$.

The construction will be different according to the congruence class of $\ell$  modulo 8.

\begin{lemma}\label{le:C01}
Let $\ell \ge  7$.  There exist two $\ell$-cycles, $C_0$ and $C_1$, 
 on vertex set $\Z_{2\ell}\times\{0, 1\}$ having among them
 the two $2\ell-2$ edges of mixed differences in $\{0,1,\dots,2\ell-1\}\setminus\{0,\ell\}$ and one edge of $i$-pure difference $\pm m$ in $C_i$, $i=0,1$, where $m$ is one of $(\ell-5)/2,(\ell-3)/2,(\ell-1)/2,(\ell+1)/2$.  Moreover $C_i$ ($i=0,1$) and its translates when developed modulo $(2\ell,-)$ are  all equitably $2$-colourable under the colouring previously described.
\end{lemma}
\begin{proof}

We will use the following notation: for $x \in \mathbb{Z}_{4\ell}$, $i \in \mathbb{Z}_2$ and $a \in \{0, 1, \ldots, 4\ell-1\}$, let $p_i(x;a)$ be the 2-path $x_i, (x+a)_{i+1}, (x-2)_i$, with subscripts computed modulo 2. For instance, $p_0(5,4)=[5_0,9_1,3_0]$ 

Let $z_i(x;a,b)$ denote the concatenation of $b$ 2-paths $p_i(x;a)$, $p_i(x-2;a+4)$, $p_i(x-4,a+8)$, $\ldots$, $p_i(x-2b+2;a+4b-4)$ (if $b=0$ this is just an empty path).  Note that $z_i(x;a,b)$ is a path of length $2b$ with vertex set $(\{x-2b,x-2b+2, \ldots, x\} \times \{i\}) \cup (\{x+a, x+a+2, x+a+4, \ldots, x+a+2b-2\} \times\{i+1\})$.  Its initial vertex is $x_i$ and terminal vertex is $(x-2b)_i$.
Moreover, $z_0(x;a,b)$ and $z_1(x;a,b)$ between them contain an edge of each mixed difference in $\pm \{a,a+2, a+4, \ldots, a+4b-2\}$.   For instance, $z_0(5,4,2)=[5_0,9_1,3_0,11_1,1_0]$.

\subsubsection*{Case $\ell\equiv 7\pmod{8}$}

$C_0$ is the concatenation of paths 
\[
\begin{array}{l}
z_1((\ell-3)/2; 3, (\ell-3)/4),\\
0_1, ((\ell-3)/2)_0,\\
z_0((\ell-3)/2;3, (\ell-3)/4), \\
0_0,((\ell-1)/2)_0,((\ell-3)/2)_{1}
\end{array}
\]
Note that, except for vertex $((\ell-1)/2)_0$,  $C_0$ either contains both of $x_0$ and $x_1$, or else contains neither, so that $C_0$ and all its translates are equicoloured. Also,
$C_0$ contains an edge of $0$-pure difference $\pm(\ell-1)/2$, and an edge of each mixed difference in 
\[
\{-1\} \cup \pm \{3, 5, 7, \ldots, (\ell-2)\} \cup \{-(\ell-3)/2\} .
\]
$C_1$ is the concatenation of paths 
\[
\begin{array}{l}
z_0((\ell-3)/2; 2, (\ell-7)/8), \\
 z_0((\ell+1)/4; (\ell+1)/2, (\ell+1)/8)), \\
 0_0, ((\ell-3)/2)_1, \\
z_1((\ell-3)/2; 2, (\ell-7)/8),  \\
z_1((\ell+1)/4; (\ell+1)/2, (\ell+1)/8)), \\
0_1, ((\ell-1)/2)_1, ((\ell-3)/2)_0
\end{array}
\]
As above, $C_1$ and all its translates are equitably coloured. Also,
$C_1$ contains an edge of $1$-pure difference $\pm(\ell-1)/2$, and an edge of each mixed difference in 
\[
 \{1\} \cup \pm \{2,4, 6, \ldots,\xcancel{(\ell-3)/2},\ldots, (\ell-1)\} \cup \{ (\ell-3)/2\}.
\]

\subsubsection*{Case $\ell\equiv 1\pmod{8}$}

$C_0$ is the concatenation of paths \[
\begin{array}{l}
z_0((\ell-5)/2;3, (\ell-5)/4), \\
0_0, ((\ell-5)/2)_1,\\
z_1((\ell-5)/2;3, (\ell-5)/4), \\
0_1,
(\ell-1)_0,
((\ell-1)+(\ell-3)/2)_0,\\
\hspace*{1cm}((\ell-1)+(\ell-5)/2)_1, ((\ell-5)/2)_0\\
\end{array}
\]

Note that the three vertices $(\ell-1)_0,
((\ell-1)+(\ell-3)/2)_0$ and
$((\ell-1)+(\ell-5)/2)_1$ cannot all have the same colour, since $((\ell-1)+(\ell-3)/2)_0$ and
$((\ell-1)+(\ell-5)/2)_1$ have different colours; the same holds in all but two of the translates of this 3-vertex-set, and it is readily seen that in these two cases $(\ell-1)_0$ has a different colour.
As for the remaining $\ell-3$ vertices of the cycle say $x_i$, we have that $C_0$ either contains both of $x_0$ and $x_1$, or else contains neither, so that $C_0$ and its translates are equitably coloured. Moreover, 
$C_0$ contains an edge of $0$-pure difference $\pm(\ell-3)/2$, and an edge of each mixed difference in 
\[
\{-1, (\ell-5)/2)\} \cup \pm \{3, 5, 7, \ldots, (\ell-4)\} \cup \pm \{(\ell-1)\}.
\]

$C_1$ is the concatenation of paths 
\[
\begin{array}{l}
z_1((\ell-5)/2; 2, (\ell-9)/8), \\
 z_1((\ell-1)/4; (\ell-1)/2, (\ell-1)/8)), \\
0_1, ((\ell-5)/2)_0,\\
z_0((\ell-5)/2; 2, (\ell-9)/8),  \\
z_0((\ell-1)/4; (\ell-1)/2, (\ell-1)/8)),\\ 
0_0, (\ell-2)_1,
((\ell-2)+(\ell-3)/2)_1,\\
\hspace*{1cm}((\ell-2)+(\ell-5)/2)_0,((\ell-5)/2)_1
\end{array}
\]

Reasoning as for  $C_0$, the cycle $C_1$ is also equitably coloured: two different colours appear in the last three vertices, while for the remaining ones the cycle either contains both of $x_0$ and $x_1$, or else contains neither. 

The cycle $C_1$  contains an edge of $1$-pure difference $\pm(\ell-3)/2$, and an edge of each mixed 
difference in 
\[
\{1, -(\ell-5)/2)\} \cup \pm \{2, 4, 6, \ldots,\xcancel{(\ell-5)/2},\dots (\ell-3)\} \cup \pm \{(\ell-2)\}.
\]

\subsubsection*{Case $\ell\equiv 5\pmod{8}$}

$C_0$ is the concatenation of paths 
\[
\begin{array}{l}
z_0((\ell-5)/2;1, (\ell-5)/4),\\
0_0, ((\ell-5)/2)_1, \\
z_1((\ell-5)/2;1, (\ell-5)/4), \\
0_1,(\ell-2)_0,
((\ell-2)+(\ell-1)/2)_0,\\
\hspace*{1cm}((\ell-2)+(\ell-5)/2)_1, ((\ell-5)/2)_0 \\
\end{array}
\]
Reasoning as previously, it is easy to see that $C_0$ is equitably coloured.
It contains an edge of $0$-pure difference $\pm(\ell-1)/2$, and an edge of each mixed difference in 
\[
\{-2, (\ell-5)/2\} \cup \pm\{1,3, 5, 7, \ldots, (\ell-6)\} \cup  \pm \{(\ell-2)\}.
\]

$C_1$ is the concatenation of paths 
\[
\begin{array}{l}
z_1((\ell-5)/2; 4, (\ell-13)/8), \\
 z_1((\ell+3)/4; (\ell-1)/2, (\ell+3)/8)), \\
 0_1, ((\ell-5)/2)_0,\\
z_0((\ell-5)/2; 4, (\ell-13)/8),  \\
z_0((\ell+3)/4; (\ell-1)/2, (\ell+3)/8)),\\
0_0, 
(\ell-4)_1,
((\ell-4)+(\ell-1)/2)_1,\\
\hspace*{1cm}((\ell-4)+(\ell-5)/2)_0, ((\ell-5)/2)_1\\
\end{array}
\]

Reasoning as above, it is easy to see that also $C_1$ is equitably coloured. It contains an edge of $1$-pure difference $\pm(\ell-1)/2$,  and an edge of each mixed difference in 
\[
\{2,-(\ell-5)/2\} \cup \pm \{4, 6, \ldots,\xcancel{(\ell-5)/2},\ldots, (\ell-1)\} \cup \pm \{(\ell-4)\}.
\]

\subsubsection*{Case $\ell\equiv 3\pmod{8}$}

$C_0$ is the concatenation of paths
\[
\begin{array}{l}
z_0((\ell-7)/2;2,(\ell-11)/8), \\
((\ell-3)/4)_0, ((3\ell-13)/4)_1, ((\ell-7)/4)_0, ((3\ell-9)/4)_1, ((\ell-11)/4)_0, \\
z_0((\ell-11)/4; (\ell+13)/2, (\ell-11)/8), \\
0_0, ((\ell-7)/2)_1, \\
z_1((\ell-7)/2;2,(\ell-11)/8), \\
((\ell-3)/4)_1, ((3\ell-13)/4)_0, ((\ell-7)/4)_1, ((3\ell-9)/4)_0, ((\ell-11)/4)_1, \\
z_1((\ell-11)/4; (\ell+13)/2, (\ell-11)/8), \\
0_1, (2\ell-1)_0, ((\ell-7)/2)_0
\end{array}
\]
$C_0$ contains an edge of $0$-pure difference $\pm(\ell-5)/2$, and an edge of each mixed difference in
\[
\begin{array}{l}
\{1, (\ell-7)/2 \} \cup \pm \{2, 4, \ldots, (\ell-11)/2\}
\\
\hspace*{1cm} \cup \pm \{ (\ell-5)/2, (\ell-3)/2, (\ell-1)/2, (\ell+1)/2\} 
\\
\hspace*{1cm} \cup \pm \{ (\ell+13)/2, (\ell+17)/2, \ldots, \ell-1\}.
\end{array}
\]

$C_1$ is the concatenation of paths
\[
\begin{array}{l}
z_1((\ell-7)/2;3, (\ell-11)/8), \\
((\ell-3)/4)_1, ((3\ell+3)/4)_0, ((\ell-7)/4)_1, ((3\ell+7)/4)_0, ((\ell-11)/4)_1, \\
z_1((\ell-11)/4; (\ell+11)/2, (\ell-11)/8), \\
0_1, ((\ell-7)/2)_0, \\
z_0((\ell-7)/2;3, (\ell-11)/8), \\
((\ell-3)/4)_0, ((3\ell+3)/4)_1, ((\ell-7)/4)_0, ((3\ell+7)/4)_1, ((\ell-11)/4)_0, \\
z_0((\ell-11)/4; (\ell+11)/2, (\ell-11)/8), \\
0_0, (2\ell-1)_1, ((\ell-7)/2)_1
\end{array}
\]

$C_1$ contains an edge of $1$-pure difference $\pm(\ell-5)/2$, and an edge of each mixed difference in
\[
\begin{array}{l}
\{-1, -(\ell-7)/2 \} \cup \pm \{3, 5, \ldots, (\ell-9)/2\}
\\
\hspace*{1cm} \cup \pm \{ (\ell+3)/2, (\ell+5)/2, (\ell+7)/2, (\ell+9)/2\} 
\\
\hspace*{1cm} \cup \pm \{ (\ell+11)/2, (\ell+15)/2, \ldots, \ell-2\}.
\end{array}
\]

\end{proof}

\begin{example}
Case $\ell\equiv 7 \pmod{8}$: for $\ell=15$ the cycles are 
\[
\begin{array}{ll}
C_0 &= (\blue{6_1},\red{9_0},\blue{4_1},\red{11_0},\blue{2_1},\red{13_0},\blue{0_1},\red{6_0},\blue{9_1},\red{4_0},\blue{11_1},\red{2_0},\blue{13_1},\red{0_0},\red{7_0})\\
C_1 &= (\red{6_0},\blue{8_1},\red{4_0},\blue{12_1},\red{2_0},\blue{14_1},\red{0_0},\blue{6_1},\red{8_0},\blue{4_1},\red{12_0},\blue{2_1},\red{14_0},\blue{0_1},\blue{7_1})\\
\end{array}
\]  

\noindent
An example for the case $\ell\equiv 1\pmod{8}$ is part of Example \ref{ex:general4l+1}.

\smallskip

\noindent
Case $\ell\equiv 5 \pmod{8}$: for $\ell=21$ the cycles are
  \[
\begin{array}{ll}
C_0 &= (\red{8_0},\blue{9_1},\red{6_0},\blue{11_1},\red{4_0},\blue{13_1},\red{2_0},\blue{15_1},\red{0_0}, \blue{8_1},\red{9_0},\blue{6_1},\red{11_0},\blue{4_1},\red{13_0},\blue{2_1},\red{15_0},\blue{0_1},\red{19_0},\blue{29_0},\red{27_1})\\
C_1 &= (\blue{8_1},\red{12_0},\blue{6_1},\red{16_0},\blue{4_1},\red{18_0},\blue{2_1},\red{20_0},\blue{0_1}, \red{8_0},\blue{12_1},\red{6_0},\blue{16_1},\red{4_0},\blue{18_1},\red{2_0},\blue{20_1},\red{0_0},\blue{17_1},\red{27_1},\blue{25_0})\\
\end{array}
\]  

\noindent
Case $\ell \equiv 3 \pmod{8}$: for $\ell=19$ the cycles are 
\[
\begin{array}{ll}
C_0 &= (\red{6_0}, \blue{8_1}, \red{4_0}, \blue{11_1}, \red{3_0}, \blue{12_1}, \red{2_0}, \blue{18_1}, \red{0_0}, \blue{6_1}, \red{8_0}, \blue{4_1}, \red{11_0}, \blue{3_1}, \red{12_0}, \blue{2_1}, \red{18_0}, \blue{0_1}, \blue{37_0}) \\
C_1 &= (\blue{6_1}, \red{9_0}, \blue{4_1}, \red{15_0}, \blue{3_1}, \red{16_0}, \blue{2_1}, \red{17_0}, \blue{0_1}, \red{6_0}, \blue{9_1}, \red{4_0}, \blue{15_1}, \red{3_0}, \blue{16_1}, \red{2_0}, \blue{17_1}, \red{0_0}, \red{37_1})
\end{array}
\]

\end{example}

The following results are used to construct the base cycle $C_p$ in the proof of Theorem \ref{th:4l+1} for $\ell \geq 9$.

\begin{lemma}\label{le:Cp}
Let $\ell \geq 7$ be odd. 
If $\ell \geq 9$, then for any $m \in \{ (\ell-3)/2, \ldots, \ell-3\}$, there is an $\ell$-cycle $C$ on vertex set $\mathbb{Z}_{2\ell} \times \{0,1\}$ containing $(\ell-3)/2$ edges of distinct $0$-pure difference, $(\ell-1)/2$ edges of distinct $1$-pure difference and two edges of mixed differences $0$ and $\ell$, satisfying the following properties:
\begin{enumerate}
\item $C$ contains no edges of $0$-pure difference $\pm 2$, $\pm m$ or $\ell$, and no edge of $1$-pure difference $\pm m$ or $\ell$. 
\item For every vertex $x_0$ in $C$, $x_1$ is also a vertex in $C$.  In particular, $C$ and its translates when developed modulo $(2\ell,-)$ are  all equitably $2$-colourable under the colouring previously described.
\end{enumerate}
If $\ell=7$, then there is a cycle satisfying properties 1 and 2 with $m=3$.  
\end{lemma}

\begin{proof}
In the case $\ell=7$, we take
\[
C=({\color{red}0_0}, {\color{red}1_0}, {\color{red}5_0}, {\color{blue}5_1}, {\color{blue}0_1}, {\color{blue}1_1}, {\color{red}7_1}),
\]
and note that $C$ contains no edge of $0$- or $1$-pure difference $m=\pm3 = \pm(\ell-1)/2$.

Now suppose $\ell \geq 9$.  Notationally, let ${y}_i(x;a,b)$ be the path
\[
x_i, (x+a+b)_i, (x+1)_i, (x+a+b-1)_i, \ldots, t_i,
\]
where 
\[
t= x +\left\{
\begin{array}{ll}
a+b/2, \; & \mbox{if $b$ is even,} \\
(b+1)/2, & \mbox{if $b$ is odd.}
\end{array}
\right.
\]
In particular, $y_i(x;a,b)$ contains an edge of each $i$-pure difference $a, a+1, \ldots, a+b$.
Let ${y}_i'(x;a,b)$ denote the path ${y}_i(x;a,b)$ traversed in the opposite order, that is, with initial vertex $t_i$ and terminal vertex $x_i$.

If $\ell > 9$, let $Y_0={y}_0(2;3,(\ell-11)/2)$ and $Y_1={y}_1'(2;3,(\ell-11)/2)$.  For $\ell=9$, let $Y_0$ and $Y_1$ be null paths. Note that $Y_i$ is a path containing $(\ell-9)/2$ edges of distinct $i$-pure difference in $\{3, \ldots, (\ell-5)/2\}$, and vertices distinct from $0_i$, $1_i$ and $\ell_i$.  Moreover, the vertices of $Y_1$ are precisely those of the form $v_1$ where $v_0 \in V(Y_0)$.  Let the terminal vertex $Y_0$ (resp.\ the initial vertex of $Y_1$) be $t_0$ (resp.\ $t_1$).  

Let $u$ be any element of $\{(\ell+1)/2, \ldots, \ell-2\}$ such that $m \notin \{u-2,u\}$.  Note that the number of elements in this set is $(\ell-3)/2$, and at most two of them can equal $m$ or $m+2$; thus since
$\ell \geq 9$, a suitable element $u$ can be chosen.

Define the cycle $C$ as the concatenation of the following paths:
\[
\begin{array}{l}
0_0 1_0 \ell_0 2_0; \\
Y_0 \\
t_0 t_1 \\
Y_1 \\
2_1 u_1 0_1 1_1 \ell_1 
\end{array}
\]
By construction, $C$ has the required properties.
\end{proof}

\begin{example}
Let $\ell=13$ and $m=5$.  Then, choosing $u=8$, we get
\[
C=({\color{red}0_0}, {\color{red}1_0}, {\color{blue}13_0}, {\color{red}2_0}, {\color{red}6_0}, {\color{red}3_0}, {\color{blue}3_1}, {\color{blue}6_1}, {\color{blue}2_1}, {\color{blue}8_1}, {\color{blue}0_1}, {\color{blue}1_1}, {\color{red}13_1}).
\]
\end{example}

\begin{lemma}\label{le:l=3mod8}
    Let $\ell \equiv 3$ (mod $8$), $\ell \geq 11$.  There is an $\ell$-cycle $C$ on vertex set $\mathbb{Z}_{2\ell} \times \{0,1\}$ containing $(\ell-3)/2$ edges of distinct $0$-pure difference, $(\ell-1)/2$ edges of distinct $1$-pure difference and two edges of mixed differences $0$ and $\ell$, satisfying the following properties:
\begin{enumerate}
\item $C$ contains no edges of $0$-pure difference $\pm 2$, $\pm(\ell-5)/2$ or $\ell$, and no edge of $1$-pure difference $\pm(\ell-5)/2$ or $\ell$. 
\item For every vertex $x_0$ in $C$, $x_1$ is also a vertex in $C$.  In particular, $C$ and its translates when developed modulo $(2\ell,-)$ are  all equitably $2$-colourable under the colouring previously described.
\end{enumerate}
\end{lemma}

\begin{proof}
For $\ell=11$, let 
\[
C=(\red{0_0}, \red{1_0}, \blue{11_0}, \red{2_0}, \red{6_0}, \blue{6_1}, \blue{2_1}, \blue{7_1}, \blue{0_1}, \blue{1_1}, \red{11_1}).
\]

For $\ell>11$, the cycle is constructed similarly to that in Lemma~\ref{le:Cp}, except that we skip the edge of pure difference $\pm (\ell-5)/2$.
 Define the notation $y_i(x;a,b)$ and $y_i'(x;a,b)$ as in the proof of Lemma~\ref{le:Cp}.  Let $Y_0=y_0(4;3,(\ell-15)/2)$ and $Y_1=y_1'(4;3,(\ell-15)/2)$; 
 here $b=(\ell-15)/2$ is even so that $t=4+3+(\ell-15)/4 = (\ell+13)/4$.

The cycle $C$ is the concatenation of the following paths:
\[
\begin{array}{l}
    0_0 1_0 \ell_0 2_0 ((\ell+1)/2)_0 4_0;\\
    Y_0;\\
    t_0 t_1;\\
    Y_1;\\
    4_1 ((\ell+1)/2)_1 2_1 ((\ell+3)/2)_1 0_1 1_1 \ell_1.
\end{array}
\]
By construction, $C$ has the required properties.
\end{proof}

\begin{example}
    Let $\ell=19$.  Then
    \[
    C=(\red{0_0}, \red{1_0}, \blue{19_0}, \red{2_0}, \red{10_0}, \red{4_0}, \red{9_0}, \red{5_0}, \red{8_0}, \blue{8_1}, \blue{5_1}, \blue{9_1}, \blue{4_1}, \blue{10_1}, \blue{2_1}, \blue{11_1}, \blue{0_1}, \blue{1_1}, \red{19_1}).
    \]
\end{example}

Finally, the following lemma describes how to build the cycles containing $\infty$.

\begin{lemma} \label{CyclesThroughInfinity}
Let $\ell$ be odd and let $D=\{d_1, d_2, \ldots, d_{(\ell-3)/2}\}$ be a set of $(\ell-3)/2$ distinct elements of $\{1, 2, \ldots, \ell-1\}$ with $d_1 > d_2 > \cdots > d_{(\ell-3)/2}$.  There exists an $\ell$-cycle $C$ on vertex set $\mathbb{Z}_{2\ell} \cup \{\infty\}$ with $\partial(C)= \pm D \cup \{\ell\}$.
\end{lemma}

\begin{proof}
For $k=1, 2, \ldots, (\ell-3)/2$, let $s_i = \sum_{i=1}^k (-1)^{i-1} d_i$.  Since $d_1, d_2, \ldots, d_{(\ell-3)/2}$ is a strictly decreasing sequence of positive integers and each $d_i \leq \ell-1$, note that for each $i$, $0 < s_i < \ell$.  Thus
\[
C=(\infty, 0, s_1, s_2, \ldots, s_{(\ell-3)/2}, s_{(\ell-3)/2}+\ell, s_{(\ell-5)/2}+\ell, \ldots, s_1+\ell, \ell)
\]
is an $\ell$-cycle with vertices in $\mathbb{Z}_{2\ell} \cup \{\infty\}$, and it is easy to check that $\partial(C)$ is as required.
\end{proof}
\begin{remark}
    If using the $2$-colouring of $\Z_{2\ell}$ described at the beginning of this section, that is 
    one colour on $0,\dots,\ell-1$ and a different colour on $\ell,\dots,2\ell-1$, then the cycle built in Lemma \ref{CyclesThroughInfinity} is equitably coloured.
\end{remark}

\section{Equitably 2-colourable $\ell$-cycle decompositions of $K_v$ when $v \equiv 1$ or $\ell \pmod{2\ell}$}

The aim of this section is to complete the proof of the existence of an equitably 2-colourable $\ell$-cycle decomposition of $K_{v}, v\equiv 1$ or $\ell \pmod{2\ell}$. The case $v=2\ell+1$ and $v=4\ell +1$ have been already established, so if $v \equiv 1 \pmod{2\ell}$, we can assume $v\ge 6\ell+1$. 

Throughout this section, we make use of lexicographic products.  Formally, the {\em lexicographic} product of graphs $G$ and $H$ is the graph $G[H]$ with vertex set $V(G) \times V(H)$ and edges of the form $(g,h)(g',h')$, where $gg' \in E(G)$ and $h,h' \in V(H)$ or $g=g'$ and $hh' \in E(H)$.  We will generally take $H$ to be an empty graph, i.e.\ $H\simeq\overline{K_n}$, and in this case, we denote $G[H]$ by $G[n]$; the formation of $G[n]$ can be thought of as ``blowing up'' the vertices of $G$ by a factor of $n$.

Our approach for forming an equitably 2-colourable decompositions in the case $v \equiv 1 \pmod{2\ell}$ is as follows.  Let $v=2k\ell+1$, and view the vertex set of $K_v$ as $(\mathbb{Z}_{2k} \times \mathbb{Z}_{\ell}) \cup \{\infty\}$.  We partition $\mathbb{Z}_{2k} \times \mathbb{Z}_{\ell}$ into $2k$ parts of size $\ell$. Each part will be coloured with $(\ell+1)/2$ red vertices and $(\ell-1)/2$ blue vertices, while $\infty$ will be coloured blue.  The parts are then partitioned into pairs; each pair, together with $\infty$, induces a copy of $K_{2\ell+1}$, which is decomposed into equitably coloured $\ell$-cycles.  The remaining edges induce a subgraph isomorphic to $K_{k}[2\ell]$, which we decompose into copies of $C_3[\ell]$ and $C_5[\ell]$.  These are further decomposed into equitably 2-coloured $\ell$-cycles.

When $v \equiv \ell \pmod{2\ell}$, writing $v=(2k+1)\ell$, we decompose $K_v$ into $K_{2k+1}[\ell]$ together with $2k+1$ copies of $K_{\ell}$.  Again, each part is coloured with $(\ell+1)/2$ red vertices and $(\ell-1)/2$ blue vertices.  Each copy of $K_{\ell}$ is decomposed into equitably coloured Hamiltonian cycles, while $K_{2k+1}[\ell]$ is again decomposed into copies of $C_3[\ell]$ and $C_5[\ell]$, which we further decompose into equitably 2-coloured $\ell$-cycles.

Our first aim will be to show that the graphs $C_3[\ell]$ and $C_5[\ell]$ possess an equitably $2$-colourable $\ell$-cycle system; many of the ideas and constructions used here are based on the results by Alspach, Schellenberg, Stinson and Wagner in \cite{ASSW}, and in fact, in most cases our proofs will only show that the cycle decompositions built in \cite{ASSW} are 2-equicolourable. Let us point out that, unlike the authors of \cite{ASSW}, our aim is to only build a cycle system, not a 2-factorization. 

When decomposing $C_s[\ell]$, it is useful to view it as a Cayley graph.  
\begin{definition}
Let $G$ be an additive group and let $\Omega \subseteq G$ be closed under negation.  The {\em Cayley graph} $\Cay[G,\Omega]$ is the graph with vertex set $G$, such that $gh$ is an edge if and only if $g-h \in \Omega$.  
\end{definition}

Identifying the vertex set of $C_s[\ell]$ with $\Z_s \times \Z_{\ell}$, note that 
$C_s[\ell] \simeq \Cay[\Z_s\times\Z_\ell, \Omega],$ with $\Omega=\{(\pm 1,\pm i),i=0,1,\dots,(\ell-1)/2\}$.
As in \cite{ASSW}, to obtain an $\ell$-cycle system of $C_s[\ell]$, we will partition $\Omega$ into disjoint subsets $\Omega_1$ and $\Omega_2$, and decompose the graphs $\Cay[\Z_s\times\Z_\ell, \Omega_1]$ and $\Cay[\Z_s\times\Z_\ell, \Omega_2]$ separately into $\ell$-cycles; often $\Omega_1=\{(\pm 1,\pm i) \mid i=0,1,2\}$ and 
$\Omega_2=\{(\pm 1,\pm i) \mid i=3,\dots,(\ell-1)/2\}$.

\begin{definition}
Given an $\ell$-cycle $C=((a_0,b_0),(a_1,b_1),\dots,(a_{\ell-1},b_{\ell-1}))$ in $C_s[\ell]$, the $s$-development of $C$ is the set of cycles
$$\{((a_0+i,b_0),(a_1+i,b_1),\dots,(a_{\ell-1}+i,b_{\ell-1})) \mid i \in \Z_s\}$$ that is, the set of cycles obtained from $C$ by developing modulo $(s,-)$.  
\end{definition}

To decompose $C_s[\ell]$, as in~\cite{ASSW}, we employ a Hamiltonian decomposition of an auxiliary Cayley graph $\Cay[\Z_{\ell},A]$ and lift these cycles to $\Cay[\Z_s \times \Z_{\ell}, \{\pm 1\} \times A]$.  

\begin{definition}\cite{ASSW, BDT2017}
Let $C=(c_0,c_1,\dots,c_{\ell-1})$ be a directed Hamiltonian cycle in $K_\ell$ (with vertex set identified with $\Z_\ell$). 
Given the directed Hamiltonian $\ell$-cycle $C=(c_0,c_1,\dots,c_{\ell-1})$, the {\em projection} of $C$ onto $C_s[\ell]$ 

is the following $\ell$-cycle in $C_s[\ell]$:
\[
((0,c_0),(1,c_1),\dots,(s-1,c_{s-1}),(0,c_s),(1,c_{s+1}),\dots,(0,c_{\ell-2}),(1,c_{\ell-1})).
\]
The {\em reverse projection} of $C$ is the $\ell$-cycle
\begin{center}
$((0,c_0), (s-1,c_1), (s-2,c_2),\ldots,(1,c_{s-1}),(0,c_s),(s-1,c_{s+1}),(0,c_{s+2}),$ 

\hspace*{2in}
$(1,c_{s+3}), \ldots, (s-1,c_{\ell-1}))$.
\end{center}
\end{definition}

Informally, we can think of the projection of $C$ as wrapping once around the $s$-cycle, and then zig-zagging between the first two parts, while the reverse projection wraps in the opposite direction.  An edge $xy$ of difference $d=\pm(x-y)$ in $C$ will yield edges of differences $\pm(1,d)$ in $C_s[\ell]$.

It is straightforward to check that the following lemma holds.  For the decomposition part, see Lemma 2 in \cite{ASSW} and Lemma 2.6 in~\cite{BDT2017}.  

\begin{lemma}\label{lemma: projections}

If $\cal C$ is a Hamiltonian decomposition for $\Cay[Z_\ell,A]$, then we obtain an $\ell$-cycle decomposition for  $\Cay[\Z_s\times\Z_\ell, \Omega]$ with $\Omega=\{\pm 1\} \times A$ by considering the $s$-development of the projection and reverse projection of $C$ for all $C \in \cal C$.
 
Moreover, if we consistently colour each $\ell$-part with $(\ell+1)/2$ red and $(\ell-1)/2$ blue vertices, then the cycles we obtain are $2$-equicoloured.  
\end{lemma}

The lemma below easily follows from the main theorem of \cite{BFM}. 

\begin{lemma}\label{lemma: hamiltonian}
The graph $\Cay[Z_\ell,\pm\{3,4,\dots,(\ell-1)/2\}]$, $\ell\ge 7$ odd, has a Hamiltonian decomposition.
\end{lemma}
\begin{proof} 

From \cite{BFM} we have that every 4-regular connected Cayley graph on a finite abelian group can be decomposed into two Hamiltonian cycles. When $\ell\equiv 1\pmod{4}$, we obtain a  Hamiltonian decomposition of $\Cay[Z_\ell,\pm\{3,4,\dots,\\(\ell-1)/2\}]$ by applying this result iteratively to $\Cay[Z_\ell,\pm\{3,4\}]$,$\Cay[Z_\ell,\pm\{5,6\}]$, \dots,  $\Cay[Z_\ell,\pm\{(\ell-3)/2,(\ell-1)/2\}]$. If  $\ell\equiv 3\pmod{4}$ we use the same procedure to obtain a Hamiltonian decomposition for 
$\Cay[Z_\ell,\pm\{3,4,\dots,(\ell-3)/2\}]$ and note that $\Cay[Z_\ell,\pm(\ell-1)/2]$ is a Hamiltonian cycle.
\end{proof}

As in \cite{ASSW}, the general strategy to construct  $\ell$-cycle systems of $C_s[\ell]$ is to separately decompose the graphs $\Cay[\Z_s\times\Z_\ell, \Omega_1]$ and $\Cay[\Z_s\times\Z_\ell, \Omega_2]$, with $\Omega_1 \dot\cup \, \Omega_2=\Omega$, into $\ell$-cycles; with one exception (a direct construction presented in Lemma \ref{lemma: 5-cycle7}) we will have that $\Omega_1=\{(\pm 1,\pm i),i=0,1,2\}$ and 
$\Omega_2=\{(\pm 1,\pm i),i=3,\dots,(\ell-1)/2\}$. Each decomposition will be 2-equicolourable, and will have the property that each part has $(\ell+1)/2$ red and $(\ell-1)/2$ blue vertices. 

It is easy to see that such a decomposition exists for  $\Cay[\Z_s\times\Z_\ell, \Omega_2]$, since it is enough to combine Lemmas \ref{lemma: projections} and \ref{lemma: hamiltonian}: the real work will lie in decomposing $\Cay[\Z_s\times\Z_\ell, \Omega_1]$.

\begin{lemma} \label{lemma: 3-cycle}
There is an equitably $2$-colourable $\ell$-cycle decomposition of $C_3[\ell]$ in which each part has $(\ell+1)/2$ red and $(\ell-1)/2$ blue vertices.
\end{lemma}

\begin{proof}

As mentioned above, we only need to decompose $\Cay[\Z_3\times\Z_\ell, \Omega_1]$, with $\Omega_1=\{(\pm 1,\pm i),i=0,1,2\}$,  into $\ell$-cycles, and show that such a decomposition is equitably $2$-colourable with each part having $(\ell+1)/2$ red and $(\ell-1)/2$ blue vertices.

We obtain the required decomposition by developing the following five base cycles modulo 3 with respect to the first components, where, in the last cycle, $a\equiv (\ell-1)/2 \pmod{3}$.
\begin{align*}
C_1=&((0,0),(1,\ell-1),(0,\ell-2),(1,\ell-3),\dots,(1,4),(0,3),(1,2),(2,2)),\\
C_2=&((0,1),(1,0),(2,0),(1,\ell-2),(2,\ell-2),(1,\ell-4),(2,\ell-4),\dots,(1,3),(2,3)),\\
C_3=&((0,2),(1,1),(2,1),(0,\ell-1),(2,\ell-1),(0,\ell-3),(2,\ell-3),\dots,(0,4),(2,4)),\\
C_4=&((0,0),(1,1),(2,2),(0,3),(2,4),(0,5),(2,6),\dots,(0,\ell-2),(2,\ell-1)),\\
C_5=&((0,0),(1,2),(2,4),(0,6),(1,8),\dots,(a,\ell-1),\\
         &(a+1,1),(a+2,3),(a+1,5),(a,7),\dots,(0,\ell-6,),(2,\ell-4),(1,\ell-2)).	
\end{align*}
Note that this is just a different presentation of the construction used in the proof of Theorem 5 in \cite{ASSW}. We will colour each $\ell$-part of the graph
$C_3[\ell]$ with the same pattern: to ensure that the colouring is equitable, then, we only need to consider the second components of the vertices in each base cycle.
When $\ell\equiv 1 \pmod{4}$, we will colour the vertices  in $\{(i,0),(i,1),\dots,(i,(\ell-1)/2), i\in\Z_3\}$ red and the remaining vertices blue; when $\ell\equiv 3 \pmod{4}$, we will instead colour  $\{(i,0),(i,1),\dots,(i,(\ell-3)/2), i\in\Z_3\}$ blue and the remaining vertices red; it is easy to check that this colouring is equitable.

\end{proof}

When considering $\ell$-cycle decompositions for $C_5[\ell]$, the case $\ell=7$ needs an ad hoc construction, while for $\ell>7$ we once more rely on a construction from \cite{ASSW}.
\begin{lemma} \label{lemma: 5-cycle7}
There is an equitably $2$-colourable $7$-cycle decomposition of $C_5[7]$ in which each part has four red vertices and three blue vertices.
\end{lemma}

\begin{proof}

In this case, we will give an explicit cycle decomposition: we colour the vertices in the set $\{(i,2j),i\in \Z_5, j=0,1,2,3\}$ red, and the remaining vertices blue. Decompose $\Cay[\Z_5\times\Z_\ell, \{(\pm1,\pm2)\}]$	by considering the $5$-development of the projection and reverse projection of the cycle $C=(0,2,4,6,1,3,5)$, and 	$\Cay[\Z_5\times\Z_\ell, \{(\pm1,0), (\pm1,\pm1), (\pm1,\pm3)\}]$ with the $5$-development of the following five starter cycles.

\begin{align*}
C_1&=((0,1),(1,2),(2,2),(3,3),(4,0),(0,0),(1,1)),\\
C_2&=((0,3),(1,2),(2,5),(3,2),(4,6),(0,6),(1,3)),\\
C_3&=((0,5),(1,6),(2,0),(3,3),(4,4),(0,4),(1,5)),\\
C_4&=((0,1),(1,0),(2,4),(3,3),(4,6),(0,5),(1,4)),\\
C_5&=((0,4),(1,0),(2,6),(3,2),(4,1),(0,5),(1,1)).	
\end{align*}

\end{proof}

\begin{lemma}\label{lemma: 5-cycle}
There is an equitably $2$-colourable $\ell$-cycle decomposition of $C_5[\ell]$ for any odd $\ell\ge7$ in which each part has $(\ell+1)/2$ red and $(\ell-1)/2$ blue vertices.
\end{lemma}

\begin{proof}
The case $\ell=7$ follows from the previous lemma; for $\ell>7$, we consider once more the
construction used in the proof of Theorem 5 in \cite{ASSW}.
This construction, though lacking the compact presentation seen for $C_3[\ell]$, has the same pattern of that case in the second components of the vertices of the cycles.  It is then possible to see that the same colouring presented in the proof of Lemma \ref{lemma: 3-cycle} will once more give an equitably $2$-colourable cycle system.
\end{proof}

We now turn to constructing equitably 2-colourable $\ell$-cycle decompositions of $K_{v}$.

\begin{theorem} \label{th:1mod2l}
If $v \equiv 1$ (mod $2\ell$), 
then there is an equitably $2$-colourable $\ell$-cycle decomposition of $K_v$.
\end{theorem}

\begin{proof}

The result is true when $v=2\ell+1$ and $v=4\ell+1$ by Theorems~\ref{th:2l+1} and~\ref{th:4l+1}, so we henceforth assume $v \geq 6\ell+1$.

Let $v=2k\ell+1$, where $k \geq 3$.  Let the vertex set of $K_v$ be $(\mathbb{Z}_{2k} \times \mathbb{Z}_{\ell}) \cup \{\infty\}$.  For each element $x \in \mathbb{Z}_{2k}$, we colour $(\ell+1)/2$ vertices of $\{x\} \times \mathbb{Z}_{\ell}$ red and $(\ell-1)/2$ blue.  Vertex $\infty$ is coloured blue.

To find an $\ell$-cycle decomposition of $K_v$, we first decompose $K_v$ into $k$ subgraphs isomorphic to $K_{2\ell+1}$ (on vertex sets $(\{2i,2i+1\} \times \mathbb{Z}_{\ell}) \cup \{\infty\}$ for each $i \in \{0, \ldots, k-1\}$), along with a copy of $(K_{2k}-I)[\ell]$ (where the 1-factor $I$ is given by the set of edges $\{\{0,1\},\{2,3\}, \ldots, \{2k-2,2k-1\}\}$).  
Each copy of $K_{2\ell+1}$ has $\ell$ blue vertices and $\ell+1$ red vertices.  On each, we place a copy of an $\ell$-cycle decomposition of $K_{2\ell+1}$ from Theorem~\ref{th:2l+1}, with the colour classes corresponding to the existing colouring.  

It remains to decompose $(K_{2k}-I)[\ell]$.  To to this, first note that since $k \geq 3$, we can write $2k(2k-2)/2=2k(k-1)=3m+5n$ for some non-negative integers $m$ and $n$; thus $K_{2k}-I$ can be decomposed into $m$ cycles of length $3$ and $n$ cycles of length $5$ by Theorem~\ref{BHPTheorem}.  Blowing up each vertex by $\ell$, we obtain a decomposition of $(K_{2k}-I)[\ell]$ into subgraphs isomorphic to $C_3[\ell]$ and $C_5[\ell]$, where each part of size $\ell$ has $(\ell-1)/2$ blue vertices and $(\ell+1)/2$ red vertices.  On each of these subgraphs, place an equitably 2-coloured $\ell$-cycle decomposition using Lemma~\ref{lemma: 3-cycle} or Lemma~\ref{lemma: 5-cycle}.
\vspace*{2ex}

\end{proof}

\begin{theorem} \label{th:lmod2l}
Let $v \equiv \ell$ (mod $2\ell$).  There is an equitably 2-colourable $\ell$-cycle decomposition of $K_v$.
\end{theorem}

\begin{proof}
Let $v=2k\ell+\ell = \ell (2k+1)$.  We form an $\ell$-cycle decomposition of $K_v$ as follows.  
First, writing $\binom{2k+1}{2}=k(2k+1)=3m+5n$ for some non-negative integers $m$ and $n$, we may decompose $K_{2k+1}$ into $m$ cycles of length $3$ and $n$ cycles of length $5$ by Theorem~\ref{BHPTheorem}.  Now blow up each vertex by a factor of $\ell$, so that that each vertex becomes a part of size $\ell$ and each cycle becomes a copy of $C_3[\ell]$ or $C_5[\ell]$.  
Colour each part with $(\ell+1)/2$ vertices red and $(\ell-1)/2$ blue.  On each copy of $C_3[\ell]$ or $C_5[\ell]$, place an equitably $2$-coloured $\ell$-cycle decomposition by Lemma~\ref{lemma: 3-cycle} or Lemma~\ref{lemma: 5-cycle}.  Finally, on each part of size $\ell$, place a Hamiltonian decomposition of $K_{\ell}$, which is clearly equitably coloured.
\end{proof}

\section{Conclusion}

In this paper, we have shown that there is an equitably $2$-colourable $\ell$-cycle decomposition of $K_v$ whenever $\ell \equiv 1$ or $\ell \pmod{2\ell}$.  A comparable result for even orders, namely that there is an equitably $2$-colourable $\ell$-cycle decomposition of $K_v-I$ for all $\ell \equiv 0$ or $2 \pmod{\ell}$, was previously shown by the authors in~\cite{BM2021} (along with various other existence results for even cycle length).  While these results do not entirely settle the spectrum of equitably $2$-colourable cycle decompositions, they do provide further evidence that equitably colourable cycle systems exist widely, in stark contrast with equitably colourable balanced incomplete block designs (see~\cite{LP}).  Based on this evidence, together with the results of~\cite{ABW2007} for small cycle lengths, we make the following conjectures.
\begin{conjecture}
    If $v \geq \ell \geq 5$ are odd integers, then there exists an equitably $2$-colourable $\ell$-cycle decomposition of $K_v$ if and only if $\ell \mid \frac{v(v-1)}{2}$.
\end{conjecture}

\begin{conjecture}
    If $v \geq \ell \geq 4$ are even integers, then there exists an equitably $2$-colourable $\ell$-cycle decomposition of $K_v-I$ if and only if $\ell \mid \frac{v(v-2)}{2}$.
\end{conjecture}

If $\ell$ is even and $v$ is odd, it is known that no equitably $2$-colourable $\ell$-cycle decomposition of $K_v$ exists~\cite{ABW2007}.  On the other hand, the existence of equitably $2$-colourable $\ell$-cycle decompositions of $K_v-I$ with $\ell$ odd and $v$ even cannot be easily ruled out; indeed, if $\ell=5$, an equitably $2$-colourable cycle decomposition of $K_v-I$ exists for every admissible even order $v$~\cite{ABW2007}.  We leave the more general existence question of equitably $2$-colourable odd cycle decomposition of the cocktail party graph as an open problem.  

Likewise, this paper has not considered equitable $c$-colourings with $c>2$.  The case $c=3$ was considered for small cycle lengths in~\cite{ABLW2004}; nevertheless, the question of when there exists an equitably $c$-colourable $\ell$-cycle decomposition of $K_v$ or $K_v-I$ remains wide open in general.

\section*{Acknowledgements}
We thank the referees for their careful reading of the paper, for many helpful comments, and for bringing references~\cite{WangCao, WangLuCao} to our attention.

A.\ Burgess gratefully acknowledges support from NSERC Discovery Grant RGPIN-2019-04328. 
F. Merola gratefully acknowledges support from INdAM – GNSAGA Project, CUP\_E55F22000270001.

\end{document}